\documentclass{icmart}

\contact[tamarz@math.huji.ac.il]{Einstein Institute of Mathematics, Edmond J. Safra Campus, Givat Ram The Hebrew University of Jerusalem, Jerusalem, 91904, Israel; Mathematics Department, Technion - Israel Institute of Technology Haifa, 32000, Israel.}

\usepackage{tikz}
\usetikzlibrary{matrix,arrows,shapes,calc,intersections}

\usepackage{amssymb}
\usepackage{amscd}
\usepackage{amsfonts}
\usepackage{setspace}
\usepackage{version}

\newtheorem{theorem}{Theorem}[section]

\newtheorem{proposition}[theorem]{Proposition}
\newtheorem{conjecture}[theorem]{Conjecture}

\theoremstyle{definition}
\newtheorem{definition}[theorem]{Definition}
\newtheorem{remark}[theorem]{Remark}

\newcommand\E{\mathbb{E}}
\newcommand\Z{\mathbb{Z}}

\newcommand\C{\mathbb{C}}
\newcommand\N{\mathbb{N}}

\newcommand\X{\mathrm{X}}

\newcommand\D{\mathcal{D}}
\newcommand\B{\mathcal{B}}

\newcommand\F{\mathbb{F}}

\newcommand\GI{\operatorname{GI}}

\newcommand\eps{\varepsilon}

\title[Linear equations in primes and dynamics of nilmanifolds]{Linear equations in primes and dynamics of nilmanifolds}

\author[Tamar Ziegler]
{Tamar Ziegler \thanks{The author is supported by ISF grant 407/12. %
}
}

\begin{document}

\begin{abstract}
In this paper we survey some of the ideas behind the recent developments in additive number theory, combinatorics and ergodic theory leading to the proof of Hardy-Littlewood type estimates for the number of prime solutions to systems of linear equations of finite complexity.
\end{abstract}



\maketitle

\section{Introduction }

A famous conjecture of Hardy and Littlewood \cite{HL} predicts  that given a $k$-tuple of integers $\mathcal H = \{h_1, \ldots, h_k\}$, there are infinitely many $k$-tuples
\[
    x+ h_1, \ldots, x + h_k, 
\]
such that all elements are simultaneously prime unless there is an obvious divisibility obstruction. Denote by $\nu_{\mathcal H}(p)$ the number of congruence classes modulo $p$ that $\mathcal H$ occupies, and call a $k$-tuple of integers admissible if $\nu_{\mathcal H}(p)<p$ for all primes $p$. Then the Hardy-Littlewood conjecture amounts to the statement that $x + h_1, \ldots, x + h_k$ are simultaneously prime infinitely often if and only if $\mathcal H$ is admissible. Moreover, they conjectured a precise formula for the asymptotic number of $k$-tuples for an admissible $\mathcal H$: Let $\mathbb P$ denote the set of primes, then 
\[
\left| \{ x \in [1,N],  \{x+h_1, \ldots, x+h_k\} \subset \mathbb P \}\right| \sim \mathfrak{S}(\mathcal H) \frac{N}{(\log N)^k}.
\]
The constant $\mathfrak{S}(\mathcal H)$ is an Euler product and is called the singular series.\footnote{ The singular series $\mathfrak{S}(\mathcal H)$ is given by the Euler product
\[
\mathfrak{S}(\mathcal H) = \prod_p \left(1-\frac{\nu_{\mathcal H}(p)}{p}\right) \left(1-\frac{1}{p}\right)^{-k}.
\]
We refer the reader to \cite{sound} for an excellent exposition of the heuristics leading to the conjecture above. We write $a(N)\sim b(N)$ if $a(N)=b(N) (1+o(1))$. 
}
While there have recently been extraordinary developments towards our understanding of gaps between primes and prime tuples \cite{gpy, zhang, maynard, polymath}, some of them presented at the current ICM, we are still far from proving this conjecture. \\

One can relax the conjecture by looking for prime points in higher rank affine sublattices of $\Z^k$. 
In a series of papers by Green-Tao \cite{gt-linear}, \cite{gt-mobius},  Green-Tao-Z \cite{gtz} we prove:
\begin{theorem}[Green-Tao-Z (2012)]\label{main}
Let $\{\psi_i(\vec x)\}_{i=1}^{k}$ be  a collection of $k$ affine linear forms in ${m}$ variables with integer coefficients,  
$\psi_i(\vec x)= \sum_{j=1}^m a_{ij}x_j +b_i$. Suppose no two forms are affinely dependent\footnote{Affine linear forms are affinely dependent is their linear parts are linearly dependent; e.g. the forms $x$ and $x+2$ are affinely dependent. A collection of $k$ affine linear forms no two forms are affinely dependent is said to be of {\em finite complexity} \cite{gt-linear}.}. Then  
\[
\left| \{ \vec x \in [0,N]^m,  \{\psi_1(\vec x), \ldots,\psi_k(\vec x) \} \subset  \mathbb P \}\right| \sim \mathfrak{S}(\vec \psi) \frac{N^m}{(\log N)^k}
\]
where $ \mathfrak{S}(\vec \psi)$ is an explicit Euler product (analogous to $\mathfrak{S}(H)$).
\end{theorem}
As a special case of this theorem we obtain the asymptotic number of $k$-term arithmetic progressions of primes. 
The reader will observe that the condition that no two forms are affinely dependent rules out the important case of twin primes, or  more generally any $k$-tuple with bounded gaps as described above, however its non-homogeneous nature allows one to use it in various applications that were previously conditional on the Hardy-Littlewood conjectures (see for example \cite{bm,hsw}).  Theorem \ref{main} may be viewed as a vast generalization of Vinogradov's $3$-prime theorem \cite{vinogradov}:  any large enough odd number is a sum of three primes. We remark that very recently Vinogradov's result has been extended to include all odd numbers greater that $5$ \cite{helfgott}, thus verifying the weak Goldbach conjecture.

In this paper we give an outline of intertwining developments in ergodic theory, combinatorics and additive number theory leading to Theorem \ref{main}. 

\section{Arithmetic progressions in sets of positive density}
Our starting point on the combinatorial front is the following result of K. Roth \cite{roth}. Let $E \subset \N$. The {\em upper density} of $E$ is defined to be
\[
\bar{d}(E)= \limsup_{N \to \infty} \frac{|E \cap [1,N]|}{N}.
\] 

\begin{theorem}[Roth 1953] Let $E \subset \N$ be a set of positive upper density, then $E$
contains a non trivial $3$-term arithmetic progression.
\end{theorem}

Roth's proof plays an important role in later developments - we outline the idea below.  Let $\delta>0$, and suppose $E$ has density $\delta$ in an arithmetic progression $P$ of size $N$, namely $E \subset P$ and $|E| = \delta N$. We first observe that if each element in $P$ were to be chosen independently at random to be in $E$ with probability $\delta$ then $E$ would typically contain many $3$-term progressions -  approximately $\delta^3N^2$.  In view of this, Roth's argument is based on the following:
\begin{itemize}
\item either $E$ has at least $\frac{\delta^3N^2}{2}$ $3$-term progressions, or
\item $E$ has density  at least {$\delta+c(\delta)$} on a sub-progression $Q \subset P$ of size $N^{\frac{1}{3}}$, where $c$ is a decreasing positive function. 
\end{itemize}

Our starting point is a subset $E \subset P=[1,N]$, of density $\delta$. After running the above argument at most $s=1/c(\delta)$ times we obtain a subset $E' \subset E$ which is of density (exactly) $1$ in a subprogression $P'\subset [1,N]$ of size at least $N^{\frac{1}{3^s}}$.  Namely, either at some point we have many $3$-term arithmetic progressions, or after finitely many steps we find an arithmetic progression of size $N^{\frac{1}{3^s}}$ in $E$; if $N$ is sufficiently large then 
$N^{\frac{1}{3^s}} \ge 3$.

 We remark that a more careful analysis allows one to have the density $\delta$ depend on $N$ in the form $\delta = 1/(\log \log N)^{t}$ \footnote{The state of the art in the question of $3$-term progressions is the recent result of T. Sanders stating that one can have the density as small as $\delta=1/\log N^{1-o(1)}$  \cite{sanders}.}. 

The main issue is, of course, the second step in this argument - namely, obtaining increased density on a  large subprogression. This can be achieved via discrete Fourier analysis - one considers $E$ as a subset of $\Z_N=\Z/N\Z$. Denoting $1_E$ the characteristic function of $E$, one shows that if $E$ does not contain roughly the expected number of $3$-term progressions, then the function $1_E-\delta$ has a large non trivial Fourier coefficient, namely, there exist an integer $r$ such that 
\[
\left|\frac{1}{N} \sum_{x \in \Z_N} (1_E-\delta)(x)e^{2\pi i x \frac{r}{N}}\right| \ge c(\delta).
\]

Using equidistribution properties of the sequence $\{x \frac{r}{N}\}$  mod $1 $, one finds a large subprogression $Q$ - of size $N^{\frac{1}{3}}$ - on which $x \frac{r}{N}$ is roughly constant. This in turn can be translated into an increased density of at least $\delta+c(\delta)$ on (many) translates of $Q$.  This type of argument is referred to nowadays as a {\em density increment argument}. \\

Generalizing Roth's theorem to $k$-term progressions for $k>3$ turned out to be very difficult, and was shown by Szemer\'edi in his famous theorem \cite{szemeredi}:
\begin{theorem}[Szemer\'edi 1975] Let $E$ be a set of positive upper density, then $E$
contains a non trivial $k$-term arithmetic progression.
\end{theorem}

By now there are many proofs of Szemer\'edi's theorem. In this paper we will focus on two of them:  Furstenberg's ergodic theoretic proof, which marked the beginning of the ergodic theoretic side of our story, and Gowers's proof, which pioneered the application of tools from additive combinatorics to the study of arithmetic progressions.

\section{ Furstenberg's proof of Szemer\'edi's theorem.}

Shortly after Szemer\'edi proved the theorem on arithmetic progressions in sets of positive upper density in the integers, Furstenberg gave an ergodic theoretic proof of Szemer\'edi's theorem \cite{furst}. The ideas behind this proof initiated a new field in ergodic theory, referred to as {\em ergodic Ramsey theory}, and are the foundation of all subsequent ergodic theoretic developments on which the story in our paper is based.

Furstenberg first observed that one can translate questions about patterns in subsets of positive density in the integers to return time questions for sets of positive measure in a measure preserving system. More precisely:

\begin{theorem}[Furstenberg correspondence principle]
Let $\delta>0$, and let $E \subset \N$ be a set with positive upper density\footnote{Furstenberg's correspondence principle as well as his multiple recurrence theorem hold in the more general context when one considers the {\em upper Banach density} of the set $E$, $$d^*(E)=\limsup_{N-M \to \infty} \frac{|E \cap [M,N-1]|}{N-M}.$$ We will keep to the upper density for simplicity. }. There exists a probability measure preserving system\footnote{A probability {\em measure preserving system} ${\bf X}=(X, \B, \mu,T)$ consists of a probability space $(X, \B, \mu)$ and an invertible  measurable map $T:X \to X$ with $T_*\mu=\mu$. } 
 $(X, \mathcal{B}, \mu,T)$, and a measurable set $A$ with $\mu(A)>0$, such that the following holds: if for some integers $n_1, \ldots, n_k$ 
\[
\mu(A\cap T^{-n_1}A \cap \ldots \cap  T^{-n_k}A) >0,
\] 
then
\[
\bar{d}(E \cap (E-n_1) \cap \ldots \cap  (E-n_k)) >0.
\] 
In particular, there exists an integer $x$  such that $x, x+n_1, \ldots, x+n_k \in E$. 
\end{theorem}

It follows that if we seek a $k+1$ term arithmetic progression in $E$, it suffices to show that for any probability measure preserving system $(X, \mathcal{B}, \mu, T)$, and any $A$ with $\mu(A)>0$, there is a positive integer $n$ with   $\mu(A \cap T^{-n}A \cap \ldots \cap  T^{-kn}A) >0$. Observe that the case $k=1$ is the famous Poincar\'e recurrence theorem. Indeed, Furstenberg proves the following theorem:
\begin{theorem}[Furstenberg multiple recurrence theorem]\label{multiple}
Let $(X, \mathcal{B}, \mu, T)$ be a measure preserving system, and let $A$ be with $\mu(A)>0$. Then for any $k>0$
\begin{equation}\label{multi}
\liminf_{N \to \infty } \frac{1}{N} \sum_{n \le N} \mu(A \cap T^{-n}A \cap \ldots \cap  T^{-kn}A) >0.
\end{equation}
\end{theorem}
On first impression, it might seem that in replacing the arbitrary set $E$ of positive density with an arbitrary set $A$ of positive measure, our situation is not much improved. However, in the ergodic theoretic context one might hope to prove and apply useful structure theorems. In the case at hand - the averages \eqref{multi} are studied via morphism to more structured measure preserving systems, as we will try to demonstrate below. 

We will henceforth assume that the system ${\bf X}$ is {\em  ergodic}, namely any $T$-invariant set is of measure either $0$ or $1$. Any system can be decomposed to its ergodic components, thus we lose no generality in Theorem \ref{multiple} by making this assumption. 

We first briefly discuss Furstenberg's ergodic theoretic proof of Roth's theorem on $3$-term progressions. We wish to evaluate the average
\[
\frac{1}{N} \sum_{n \le N} \mu(A \cap T^{-n}A \cap  T^{-2n}A) = \frac{1}{N} \sum_{n \le N} \int 1_A(x)1_A(T^nx)1_A(T^{2n}x)d\mu
\]
where $1_A(x)$ is the characteristic function of $A$. Furstenberg proves that there exists a measure preserving system ${\bf Z}=(Z, \B_Z, \mu_Z, T_Z) $ that is a {\em Kronecker system}\footnote{A Kronecker system ${\bf Z}=(Z, \B_Z, \mu_Z, T_Z) $ is a system where $Z$ is a compact Abelian group, $\B_Z$ the Borel $\sigma$-algebra, $\mu_Z$ the Haar measure, and $T_Z$ is a rotation $T_Z(x)=x+\alpha$ for some $\alpha \in Z$}, and a morphism\footnote{A morphism between measure preserving systems ${\bf X}, {\bf Y}$ is a measure preserving map between the corresponding measure spaces that intertwines the actions of $T_X, T_Y$.  In this case ${\bf Y}$ is called a {\em factor} of ${\bf X}$.} $\pi:{\bf X} \to {\bf Z}$ such that for any $f_0, f_1, f_2 \in L^{\infty}(X)$, 
\[
\frac{1}{N} \sum_{n \le N} \int f_0(x)f_1(T^nx)f_2(T^{2n}x)d\mu
\]
is asymptotically the same as
\[
\frac{1}{N} \sum_{n \le N} \int \pi_*f_0(z)\pi_*f_1(T_Z^nz)\pi_*f_2(T_Z^{2n}z)d\mu_Z.
\]
That is, rather than trying to evaluate the average in an arbitrary (ergodic) system, we need to evaluate it in a very special system - a compact abelian group rotation system: we are left with evaluating
\[
\lim \frac{1}{N} \sum_{n \le N} \int \pi_*1_A(z)\pi_*1_A(z+n\alpha)\pi_*1_A(z+2n \alpha)d\mu_Z.
\]
Via Fourier analysis the above limit is easily seen to equal
\[
\int \pi_*1_A(z)\pi_*1_A(z+b)\pi_*1_A(z+2b)d\mu_Z(z)d\mu_Z(b).
\]
Now the projection $\pi_*$ is a positive operator, namely  if $f \ge 0$ then $\pi_*f \ge 0$. It follows that $\pi_*1_A \ge 0$, and since $\int \pi_*1_A  d\mu_Z = \int 1_A d\mu = \mu(A) >0$, the above average is clearly positive. 
The system ${\bf Z}= {\bf Z}({\bf X})$ is called the {\em Kronecker factor} of ${\bf X}$ and satisfies the following universal property.
If $\bf{Y}$ is a Kronecker system that is a factor of ${\X}$ and  $\pi_Y:{\bf X} \to {\bf Y}$ the factor map, then $\pi_{Y}$ factors through 
$ {\bf Z}({\bf X})$ as demonstrated in the diagram below:
\begin{center}
\begin{tikzpicture}
  \matrix (m) [matrix of math nodes,row sep=3em,column sep=4em,minimum width=2em]
  {
    {\bf X} & {} \\
     {\bf Z(X)} & {\bf Y} \\};
  \path[-stealth]
    (m-1-1) edge node [left] {$\pi$} (m-2-1)
            edge node [above] {$\pi_Y$} (m-2-2)
    (m-2-1.east|-m-2-2) edge [dashed] node  [above] {$\exists$} (m-2-2)
               edge [dashed,-] (m-2-1);
\end{tikzpicture}
\end{center}

The factor ${\bf Z}({\bf X})$ is constructed via the eigenfunctions of ${\bf X}$. Let us demonstrate why a non trivial eigenfunction implies the existence of a non-trivial circle rotation factor. Let  $\psi$ be an eigenfunction of ${\bf X}$,
 $$\psi(Tx) = \lambda \psi(x).$$ 
The function $|\psi|$ is a $T$-invariant function, and by ergodicity $|\psi|$ is constant a.e. Thus we can normalize $\psi$ to take values in the unit circle.  
 Any normalized eigenfunction gives rise to a morphism  to a circle rotation system $\psi: {\bf X} \to (S^1,\mathrm{Borel},\mathrm{Haar},\cdot \lambda)$:
\begin{center}
\begin{tikzpicture}
  \matrix (m) [matrix of math nodes,row sep=3em,column sep=4em,minimum width=2em, ,ampersand replacement=\&]
  {
     {\bf X} \& {\bf X}  \\
     {\bf S^1} \& {\bf S^1} \\};
  \path[-stealth]
    (m-1-1) edge node [left] {$\psi$} (m-2-1)
            edge  node [above] {$T$} (m-1-2)
     (m-1-2) edge node [left] {$\psi$} (m-2-2)        
    (m-2-1.east|-m-2-2) edge 
            node [above] {$\cdot \lambda$} (m-2-2);
\end{tikzpicture} 
\end{center}

The factor ${\bf Z}({\bf X})$ would then be the image of the map $(\psi_i): X \to (S^1)^{\mathbb N}$ given by $x \to (\psi_i(x))$, where $\{\psi_i\}$ is the collection of normalized eigenfunctions\footnote{We implicitly assume that the system ${\bf X}$ is separable and thus has at most countably many normalized eigenfunctions.} of ${\bf X}$. \\
\ \\
If ${\bf X}$ has no non-trivial eigenfunctions, then  ${\bf Z}({\bf X})$ is trivial (a point system), and thus $\pi_* f  = \int f \ d \mu$. In this case ${\bf X}$ is called {\em weakly mixing}. We then have
\[ 
 \frac{1}{N}\sum_{n=1}^N 
\int f(x)f(T^{n}x)f(T^{2n}x)\ d\mu   \to \left(\int f \ d \mu \right)^3,
\]
and we can thus think of the points $x, T^{n} x, T^{2n}x$ as asymptotically independent on average. The content of Furstenberg's argument is then that if $x, T^{n} x, T^{2n}x$ are {\em not} asymptotically independent on average, then the obstruction lies in an Abelian group rotation factor. We remark that it is clear that an Abelian group rotation factor is an obstruction as in Abelian groups  $z+2n\alpha$ is determined by $z, z+n\alpha$.
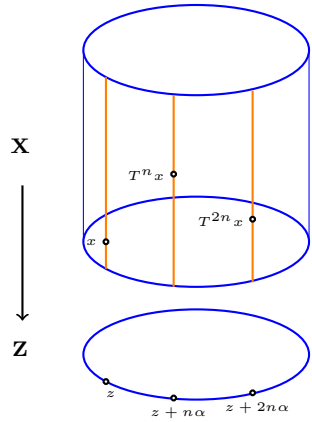
\begin{figure}[h]
\centerline{
    \begin{tikzpicture}[scale=3]
	\tikzset{mypoints/.style={fill=white,draw=black,thick}}
	\def\ptsize{0.3pt} 
	\draw[name path =ellipse,blue,thick]
      (0,0) circle[x radius = 0.5 cm, y radius = 0.2 cm];
  \draw[name path =ellipse,blue,thick]
      (0,0.5) circle[x radius = 0.5 cm, y radius = 0.2 cm];
  \draw[name path =ellipse,blue,thick]
      (0,1.35) circle[x radius = 0.5 cm, y radius = 0.2 cm];
         \coordinate (p1) at (0.5,0.5);
   \coordinate (p2) at (0.5,1.35);
   \draw[blue, thin] (p1)--(p2);
   \coordinate (p11) at (-.4,0.38);
   \coordinate (p12) at (-.4,1.23);
   \draw[orange, thick] (p11)--(p12);
   \coordinate (p21) at (-.1,0.3);
   \coordinate (p22) at (-.1,1.15);
   \draw[orange, thick] (p21)--(p22);
    \coordinate (p31) at (.25,0.32);
   \coordinate (p32) at (.25,1.17);
   \draw[orange, thick] (p31)--(p32);
      \coordinate (p01) at (-0.5,0.5);
   \coordinate (p02) at (-0.5,1.35);
   \draw[blue, thin] (p01)--(p02);
    \coordinate[label = below:{{\fontsize{8}{8} \selectfont ${\bf X}$}}](s0) at (-0.8,1);
       \coordinate[label = below:{{\fontsize{8}{8} \selectfont ${\bf Z}$}}](s1) at (-0.8,0.1);
    \coordinate (t0) at (-0.77,0.75);
   \coordinate (t1) at (-0.77,.15);
   \draw[->] [black, thick] (t0)--(t1);        
   \coordinate[label = below:{{\fontsize{4}{4} \selectfont $z$}}](p3) at (-0.4,-.12);
   \coordinate[label = below:{\fontsize{4}{4} \selectfont$z+n\alpha$}] (p4) at (-0.1,-.193);
  \coordinate[label = below:{\fontsize{4}{4} \selectfont$z+2n\alpha$}] (p5) at (0.25,-.17);
     \coordinate[label = left:{\fontsize{4}{4} \selectfont$x$}] (q3) at (-0.4,0.5);
   \coordinate[label = left:{\fontsize{4}{4} \selectfont$T^nx$}] (q4) at (-0.1,0.8);
  \coordinate[label = left:{{\fontsize{4}{4} \selectfont$T^{2n}x$}}] (q5) at (0.25,0.6);
   \foreach \p in {p3,p4,p5,q3,q4,q5}
		\fill[mypoints] (\p) circle (\ptsize);		
 \end{tikzpicture}   }
 
\caption{\small{The points $x, T^{n} x, T^{2n}x$ are independent (asymptotically on average) in the fibers  over the maximal Abelian group rotation factor. }}
\end{figure}
 
To summarize, Furstenberg's proof of Roth's Theorem is based on the following dichotomy:
\begin{itemize}
\item either ${\bf X}$ is weakly mixing, or
\item there is a morphism from ${\bf X}$ to a non trivial group rotation system. 
\end{itemize}
The above argument motivates the following definition (\cite{fw-char}):
 \begin{definition}[$k$-characteristic factor] Let ${\bf Y}$ be a factor of ${\bf X}$, and let 
$\pi: {\bf X} \to {\bf Y}$ be the factor map. We say that ${\bf Y}$ is  {\em $k$-characteristic} if 
\[ 
 \frac{1}{N}\sum_{n=1}^N 
\int f_0(x)f_1(T_{{\bf X}}^{n}x)\ldots f_k(T_{{\bf X}}^{kn}x)d\mu_{{\bf X}}  
\]
is asymptotically the same as
\[
 \frac{1}{N}\sum_{n=1}^N \int \pi_*f_0(y)\pi_*f_1(T_{{\bf Y}}^n y) \ldots \pi_*f_k(T_{{\bf Y}}^{kn}y) d\mu_{\bf Y}.
\]
\end{definition}
\ \\
We make the following observations:
\begin{itemize}
\item The system ${\bf X}$ itself is $k$-characteristic for all $k$. 
\item The trivial system is ${1}$-characteristic.  In this case $\pi_*f (x)= \int f(x)d\mu_{{\bf X}}  $, and by the {\em mean ergodic theorem}
\[
 \frac{1}{N}\sum_{n=1}^N  \int f(x)f(T_{{\bf X}}^{n}x)d\mu_{{\bf X}}  \sim \left( \int f d\mu_{{\bf X}}\right)^2.
\]
\item The Kronecker factor ${\bf Z}({\bf X})$ is ${2}$-characteristic (Furstenberg \cite{furst}).
\end{itemize}

The Furstenberg-Zimmer structure theorem \cite{furst, zimmer} relativizes the dichotomy between weak mixing and an abelian rotation factor.   One can show, using spectral theory, that ${\bf X}$ being weakly mixing is equivalent to the product system with the diagonal action ${\bf X} \times {\bf X}$ being ergodic. One can relativize this notion as follows.  Let ${\bf X} \times_{\bf Y} {\bf X}$ be the fiber product over ${\bf Y}$.  Say that $\pi:{\bf X} \to {\bf Y}$ is a {\em relatively weak mixing extension} if the map 
$\pi \times_{\bf Y}\pi: {\bf X} \times_{\bf Y} {\bf X} \to {\bf Y}$ is relatively ergodic, namely any $T \times T$ invariant subset in ${\bf X} \times_{\bf Y} {\bf X} $ is lifted  from ${\bf Y}$ via the map $\pi \times_{\bf Y}\pi $.  The role of the compact abelian group rotation is replaced by the notion of an isometric extension.
Say that $\pi:{\bf X} \to {\bf Y}$ is an {\em isometric extension} if ${\bf X}={\bf Y} \times_{\sigma} {\bf M}$ where ${\bf M}=(M, \B_{\bf M}, \mu_{\bf M})$ with $M$ a compact metric space, $\B_{\bf M}$ the Borel $\sigma$-algebra and $\mu_{M}$ the probability measure invariant under the the action of the isometry group of $M$, $T_{\bf X}(y,m)= (T_{\bf Y}y, \sigma(y)m)$, where $\sigma$ is a (measurable) map from $Y$ to the isometry group of $M$, and $\mu_{\bf X} = \mu_{\bf Y} \times \mu_{M}$.

\begin{theorem}[Furstenberg-Zimmer structure theorem \cite{furst, zimmer}]
There exists a sequence of factors 
\[
{\bf X} \to  \ldots \to  {\bf Z}_k({\bf X}) \to {\bf Z}_{k-1}({\bf X}) \to  \ldots  \to  {\bf Z}_{1}({\bf X}) \to {\star} 
\] 
such that for each $k$, either  ${\bf X} \to {\bf Z}_{k}({\bf X})$ is relatively weakly mixing, or there is a morphism from ${\bf X}$ to a non trivial isometric extension of ${\bf Z}_{k}({\bf X})$.  
\end{theorem}

\begin{theorem}[Furstenberg \cite{furst}] The factors ${\bf Z}_{k}({\bf X})$ are {$(k+1)$}-characteristic.
\end{theorem}

Observe that the factor ${\bf Z}_{0}({\bf X})$ is the trivial factor and the factor ${\bf Z}_{1}({\bf X})$ is the Kronecker factor. With the above structure theorem at hand it then suffices to prove the multiple recurrence theorem for systems which are towers of isometric extensions. Furstenberg utilizes this structure to show multiple recurrence - the idea being that if  the multiple recurrence property holds for any $k$ for a system ${\bf Y}$,  and  
${\bf X}$ is an isometric extension of ${\bf Y}$, then multiple recurrence holds for ${\bf X}$ as well.

\section{ Obstructions to $4$-term progressions.}

The Kronecker factor ${\bf Z}_{1}({\bf X})={\bf Z}({\bf X})$ is also a {\em universal ${2}$-characteristic factor} : it satisfies the property that 
if ${\bf Y}$ is any ${2}$-characteristic factor and $\pi_Y:{\bf X} \to {\bf Y}$ the factor map, then the factor map $\pi_Z:{\bf X} \to {\bf Z({\bf X})}$ factors through ${\bf Y}$ as demonstrated in the diagram  below:
\begin{center}
\begin{tikzpicture}
  \matrix (m) [matrix of math nodes,row sep=3em,column sep=4em,minimum width=2em, ,ampersand replacement=\&]
  {
     {\bf X} \& {}  \\
     {\bf Y} \& {\bf Z}({\bf X}) \\};
  \path[-stealth]
    (m-1-1) edge node [left] {$\pi_Y$} (m-2-1)
            edge  node [above] {$\pi_Z$} (m-2-2)
    (m-2-1.east|-m-2-2) edge [dashed] 
            node [above] {$\exists$} (m-2-2);
\end{tikzpicture} 
\end{center}

The factors ${\bf Z}_{k}({\bf X})$ that were constructed by Furstenberg are {\em not} universal ${(k+1)}$-characteristic for $k>1$.  This raises the following natural problem: classify the universal ${(k+1)}$-characteristic factors $Z_{k}(X)$. In other words, we try to understand the exact obstructions on the points $x, T^nx, \ldots, T^{(k+1)n}x$ preventing them from moving about freely in $X$. \\

For the case $k=1$, the upshot of the discussion regarding Furstenberg's proof of Roth's theorem on $3$-term progressions in the previous section was that the only obstructions to the independence (asymptotically on average) of $x, T^nx, T^{2n}x$ come from a compact abelian group rotation factor,  associated to the non trivial eigenfunctions of ${\bf X}$.  Already in the case $k=2$ (corresponding to $4$-term progressions) we have new obstructions. Consider for example the system 
\[
 {\bf Y}=(\mathbb T \times \mathbb T,\mathrm{Borel}, \mathrm{Haar}, T_{\bf Y})
 \]
 where
 \[
 T_{\bf Y}y=T_{\bf Y}(z, w)=(z+\alpha, w+2z+\alpha),
 \]
where $\alpha$ is irrational. Iterating $S$ we obtain
\[
 T_{\bf Y}^{n}y=T_{\bf Y}^{n}(z, w) = (z+{n}\alpha, w+{2n}z+{n^2}\alpha).
\]
We now observe that 
\[
 \begin{aligned}
y = 3T_{\bf Y}^{n}y -3T_{\bf Y}^{2n}y+T_{\bf Y}^{3n}y
 \end{aligned}
 \]
 Namely, the point $y$ is determined by the three points $T_{\bf Y}^ny, T_{\bf Y}^{2n} y, T_{\bf Y}^{3n}y$.

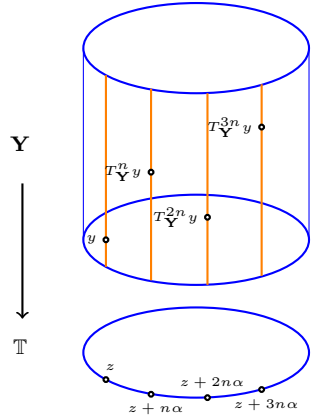
\begin{figure}[h]
\centerline{
    \begin{tikzpicture}[scale=3]
	\tikzset{mypoints/.style={fill=white,draw=black,thick}}
	\def\ptsize{0.3pt} 
      \draw[name path =ellipse,blue,thick]
      (0,0) circle[x radius = 0.5 cm, y radius = 0.2 cm];
      \draw[name path =ellipse,blue,thick]
      (0,0.5) circle[x radius = 0.5 cm, y radius = 0.2 cm];
     \draw[name path =ellipse,blue,thick]
      (0,1.35) circle[x radius = 0.5 cm, y radius = 0.2 cm];
   \coordinate (p1) at (0.5,0.5);
   \coordinate (p2) at (0.5,1.35);
   \draw[blue, thin] (p1)--(p2);
   \coordinate (p11) at (-.4,0.38);
   \coordinate (p12) at (-.4,1.23);
   \draw[orange, thick] (p11)--(p12);
   \coordinate (p21) at (-.2,0.32);
   \coordinate (p22) at (-.2,1.17);
   \draw[orange, thick] (p21)--(p22);
   \coordinate (p31) at (.05,0.3);
   \coordinate (p32) at (.05,1.15);
   \draw[orange, thick] (p31)--(p32);
   \coordinate (p31) at (.29,0.34);
   \coordinate (p32) at (.29,1.19);
   \draw[orange, thick] (p31)--(p32);
   \coordinate (p01) at (-0.5,0.5);
   \coordinate (p02) at (-0.5,1.35);
   \draw[blue, thin] (p01)--(p02);
   \coordinate[label = below:{{\fontsize{8}{8} \selectfont ${\bf Y}$}}](s0) at (-0.8,1);
   \coordinate[label = below:{{\fontsize{8}{8} \selectfont ${\bf{\mathbb T}}$}}](s1) at (-0.8,0.1);
   \coordinate (t0) at (-0.77,0.75);
   \coordinate (t1) at (-0.77,.15);
   \draw[->] [black, thick] (t0)--(t1); 
   \coordinate[label = above:{{\fontsize{4}{4} \selectfont $z$}}](p3) at (-0.4,-.12);
   \coordinate[label = below:{\fontsize{4}{4} \selectfont$z+n\alpha$}] (p4) at (-0.2,-.185);
  \coordinate[label = above:{\fontsize{4}{4} \selectfont$z+2n\alpha$}] (p5) at (0.05,-.2);
  \coordinate[label = below:{\fontsize{4}{4} \selectfont$z+3n\alpha$}] (p6) at (0.29,-.164);
  \coordinate[label = left:{\fontsize{4}{4} \selectfont$y$}] (q3) at (-0.4,0.5);
  \coordinate[label = left:{\fontsize{4}{4} \selectfont$T_{\bf Y}^ny$}] (q4) at (-0.2,0.8);
  \coordinate[label = left:{{\fontsize{4}{4} \selectfont$T_{\bf Y}^{2n}y$}}] (q5) at (0.05,0.6);
  \coordinate[label = left:{\fontsize{4}{4} \selectfont $T_{\bf Y}^{3n}y$}] (q6) at (0.29,1);
   \foreach \p in {p3,p4,p5,p6,q3,q4,q5,q6}
		\fill[mypoints] (\p) circle (\ptsize);
   \end{tikzpicture}  }
\caption{\small{The points $y, T_{\bf Y}^ny, T_{\bf Y}^{2n} y, T_{\bf Y}^{3n}y$ are not independent in the fibers over $\mathbb T$ . }}
\end{figure}

If there is a morphism $\bf{X} \to {\bf Y}$, these new obstructions to the (asymptotic on average) independence of the points $x, T_{\bf X}^nx, T_{\bf X}^{2n} x, T_{\bf X}^{3n}x$ will surface. Another way to see the obstructions coming from the system ${\bf Y}$ is by observing that the system ${\bf Y}$ exhibits {\em second order eigenfunctions}, namely functions $\phi$ satisfying  $\phi(T_{\bf Y}y)=\psi(y)\phi(y)$ where $\psi$ is an ordinary (first order) eigenfunction; for example the function $\phi(y)=\phi(z,w)=e^{2 \pi i w}$ is a second order eigenfunction. Any second order eigenfunction satisfies 
\[
\phi(y) = \phi^3(T_{\bf Y}^{n}y)\phi^{-3}(T_{\bf Y}^{2n}y)\phi(T_{\bf Y}^{3n}y)
\]
Thus choosing $f_0=\phi^{-1}$, $f_1=\phi^3$, $f_2=\phi^{-3}$, and $f_3=\phi$ we see that
\[\begin{aligned}
 1 &= \int f_0(x)f_1(T_{\bf Y}^{n}y)f_2(T_{\bf Y}^{2n}y)f_3(T_{\bf Y}^{3n}x)dm \\&= \frac{1}{N}\sum_{n \le N}\int f_0(x)f_1(T_{\bf Y}^{n}y)f_2(T_{\bf Y}^{2n}y) f_3(T_{\bf Y}^{3n}x)dm .
\end{aligned}\]
On the other hand one can verify that a (non trivial) 2nd order eigenfunction $\phi$ (and its powers) is orthogonal to ordinary eigenfunctions, thus for any $i=0,1,2,3$ the projection of the function $f_i$ on the Kronecker factor is $0$.  \\
 
It turns out, however, that second order eigenfunctions are not the only obstructions. Consider the Heisenberg nilsystem: the phase space $Y$ is the Heisenberg nilmanifold  
\[
 Y=N/\Gamma=\left( \begin{smallmatrix}
1& \mathbb R & \mathbb R  \\
0 & 1 & \mathbb R \\
0 & 0 & 1
 \end{smallmatrix} \right) \left/
 \left( \begin{smallmatrix}
1& \mathbb Z & \mathbb Z  \\
0 & 1 & \mathbb Z \\
0 & 0 & 1
 \end{smallmatrix} \right) \right. \quad  
  \]
 equipped with the Borel $\sigma$-algebra and the Haar measure, and the transformation $T_{\bf Y}$ given by $T_{\bf Y}g\Gamma = ag \Gamma$, where
 \[
 a = \left( \begin{smallmatrix}
1& \alpha & 0  \\
0 & 1 &\beta \\
0 & 0 & 1
 \end{smallmatrix} \right).
\]
Topologically $Y$ is a circle bundle over a two dimensional torus. This system shares with the system in the above example the property that  the point $g\Gamma$  is determined by $a^{{n}}g\Gamma, a^{2n}g\Gamma,a^{{3n}}g\Gamma$. However this dependence can not be described by a simple equation as in the previous example. Moreover, ${\bf Y}$ has {\em no} non-trivial second order eigenfunctions\footnote{The easiest way to see this is via equidistribution properties of polynomial orbits on nilmanifolds \cite{leibman}.}.

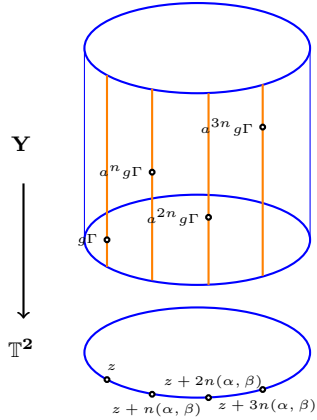
\begin{figure}[h]
\centerline{
  \begin{tikzpicture}[scale=3]
      	\tikzset{mypoints/.style={fill=white,draw=black,thick}}
	\def\ptsize{0.3pt} 
	\draw[name path =ellipse,blue,thick]
      (0,0) circle[x radius = 0.5 cm, y radius = 0.2 cm];
  \draw[name path =ellipse,blue,thick]
      (0,0.5) circle[x radius = 0.5 cm, y radius = 0.2 cm];
  \draw[name path =ellipse,blue,thick]
      (0,1.35) circle[x radius = 0.5 cm, y radius = 0.2 cm];
   \coordinate (p1) at (0.5,0.5);
   \coordinate (p2) at (0.5,1.35);
   \draw[blue, thin] (p1)--(p2);
   \coordinate (p11) at (-.4,0.38);
   \coordinate (p12) at (-.4,1.23);
   \draw[orange, thick] (p11)--(p12);
   \coordinate (p21) at (-.2,0.32);
   \coordinate (p22) at (-.2,1.17);
   \draw[orange, thick] (p21)--(p22);
    \coordinate (p31) at (.05,0.3);
   \coordinate (p32) at (.05,1.15);
   \draw[orange, thick] (p31)--(p32);
       \coordinate (p31) at (.29,0.34);
   \coordinate (p32) at (.29,1.19);
   \draw[orange, thick] (p31)--(p32);
      \coordinate (p01) at (-0.5,0.5);
   \coordinate (p02) at (-0.5,1.35);
   \draw[blue, thin] (p01)--(p02);
       \coordinate[label = below:{{\fontsize{8}{8} \selectfont ${\bf Y}$}}](s0) at (-0.8,1);
       \coordinate[label = below:{{\fontsize{8}{8} \selectfont ${\bf{\mathbb{T}^2}}$}}](s1) at (-0.8,0.1);
    \coordinate (t0) at (-0.77,0.75);
   \coordinate (t1) at (-0.77,.15);
   \draw[->] [black, thick] (t0)--(t1); 
   \coordinate[label = above:{{\fontsize{4}{4} \selectfont $z$}}](p3) at (-0.4,-.12);
   \coordinate[label = below:{\fontsize{4}{4} \selectfont$z+n(\alpha, \beta)$}] (p4) at (-0.2,-.185);
  \coordinate[label = above:{\fontsize{4}{4} \selectfont$z+2n(\alpha, \beta)$}] (p5) at (0.05,-.2);
  \coordinate[label = below:{\fontsize{4}{4} \selectfont$z+3n(\alpha, \beta)$}] (p6) at (0.29,-.164);
     \coordinate[label = left:{\fontsize{4}{4} \selectfont$g\Gamma$}] (q3) at (-0.4,0.5);
   \coordinate[label = left:{\fontsize{4}{4} \selectfont$a^ng\Gamma$}] (q4) at (-0.2,0.8);
  \coordinate[label = left:{{\fontsize{4}{4} \selectfont$a^{2n}g\Gamma$}}] (q5) at (0.05,0.6);
  \coordinate[label = left:{\fontsize{4}{4} \selectfont $a^{3n}g\Gamma$}] (q6) at (0.29,1);
   \foreach \p in {p3,p4,p5,p6,q3,q4,q5,q6}
		\fill[mypoints] (\p) circle (\ptsize);
   \end{tikzpicture}}
\caption{\small{The points $g\Gamma,  a^{{n}}g\Gamma, a^{2n}g\Gamma,a^{{3n}}g\Gamma $ are not independent in the fibers over the two dimensional torus.}}
\end{figure}

 The Heisenberg nilsystem is a special case of the following system:\\
\[{\bf Y}=(N/\Gamma, \mathrm{Borel}, \mathrm{Haar}, T_{\bf Y}),
\] where $N/\Gamma$ a {$2$}-step nilmanifold, and
\[
T_{\bf Y}: g\Gamma  \to  ag\Gamma \qquad a \in N.
\]
The system  ${\bf Y}$ is called a ${2}$-step {nilsystem}. It turns out that we need not look for further obstructions in the case $k=2$ -  all obstructions to {$4$}-term progressions come from {$2$}-step pro-nilsystems - inverse limits of $2$-step nilsystems \cite{fw-char}, \cite{cl1,cl2,cl3}:
 \begin{theorem}[Conze-Lesigne, Furstenberg-Weiss] Let ${\bf X}$ be an ergodic  measure preserving system. 
There exists a {$2$}-step pro-nilsystem ${\bf Y}$ and a morphism $\pi:{\bf X} \to {\bf Y}$ such that ${\bf Y}$ is the universal {$3$}-characteristic factor of ${\bf X}$, namely 
\[ 
 \frac{1}{N}\sum_{n=1}^N 
\int f(x)f(T_{{\bf X}}^{n}x)f(T_{{{\bf X}}}^{{2n}}x)f(T_{{{\bf X}}}^{{3n}}x)d\mu_{{\bf X}}  
\]
is asymptotically the same as
\[
 \frac{1}{N}\sum_{n=1}^N \int \pi_*f(y)\pi_*f(T_{{{\bf Y}}}^{{n}}y) \pi_*f(T_{{{\bf Y}}}^{{2n}}y)  \pi_*f(T_{{{\bf Y}}}^{{3n}}y) d\mu_{\bf Y}.
\]
\end{theorem}
We can now prove Szemer\'edi's theorem for $4$-term progressions by verifying that in a $2$-step nilsystem the above limit is positive (when $f=1_A$). \\

Let us say a few words about the proof. By Furstenberg's structure theorem it is sufficient to study systems ${\bf X}$ of the form ${\bf Z} \times_{\sigma} {\bf M}$ where ${\bf Z}={\bf Z}_1({\bf X})$ is the Kronecker factor, $T_{\bf Z}(z) = z+\alpha$, and ${M}$ is a compact metric space and $\sigma: Z \to$ ISO$({M})$. It is then shown that one can further reduce to the case where $M$ is a compact abelian group and $\sigma:Z \to M$ satisfies a functional equation now called the Conze-Lesigne equation: for all $b$, a.e $z$
\begin{equation}\label{cl}
\sigma(z+b)-\sigma(z) = c(b) + F_b(z+\alpha) - F_b(z).
\end{equation}
Describing how one can solve the above equation is beyond the scope of this paper, but let us hint how this equation is related to nilpotency. Consider the group $$G=\{(b,f): b \in Z, f:Z \to M \ \text{measurable}\}$$ with the action 
\[
(b, f) * (c, g) = (b+c, f^c\cdot g) 
\]
where $f^c(z)=f(z+c)$. 
Then condition \eqref{cl} can be interpreted as the fact that $[(\alpha, \sigma), (b, F_b)]$ is in the center of $G$, which hints at $2$-step nilpotent behavior. \\

We mention another observation regarding equation \eqref{cl}.  Upon examination one sees that 
\[
c(b_1+b_2)-c(b_1)-c(b_2)
\]
is an eigenvalue of $T_{\bf Z}$, and using the fact that there are only countably many of those, one can modify $c(b), F_b(z)$ so that $c(b)$ is linear in $b$ in a neighborhood of zero in $Z$. A similar feature will surface in the combinatorial analysis described in section \ref{inverse} below, devoted to the Inverse Theorem for the Gowers norms, which is why we mention it here. 

\section{Gowers proof of Szemer\'edi's Theorem} 

The next advancement (chronologically) was in the combinatorial front. Gowers gave a new proof for Szemer\'edi's theorem \cite{gowers}. His proof is a generalization of Roth's argument to arbitrarily long arithmetic progressions using an ingenious combination of discrete Fourier analysis and additive combinatorics; in particular Gowers obtains a Roth type bound for the density of the form $1/(\log \log N)^{c(k)}$ for some constant depending on $k$ - the length of the progression. \\

 We first fix some notation. We denote by $[N]$ the interval $[1,N]$. For a finite set $E$ we denote by $\E_{x \in E} f(x)$ the average $\frac{1}{|E|}\sum_{x \in E} f(x)$.  For two functions $f,g:[N] \to \C$ we write $f(x)\ll g(x)$ if $|f(x)| \le Cg(x)$ for some constant $C$ independent of $N$, and we write $f(x)\ll_A g(x)$ if $|f(x)| \le C(A)g(x)$ for some constant $C(A)$ independent of $N$. \\
 
 In the course of the proof Gowers defines the following norms which play a very important role in further developments.
\begin{definition}[Gowers norms] 
Let $f:\mathbb{Z}/N\mathbb{Z} \to \mathbb{C}$.  For $h \in \Z/N\Z$ define the discrete derivative in direction $h$
\[
\Delta_h f(x)=f(x+h) \overline{f(x)} 
\]
We define the $k$-th Gowers uniformity norms $U_k$ on $\mathbb C ^N$ by 
\[
\|f\|_{U^k[N]}^{2^k}=\mathbb{E}_{x, h_1, \ldots h_k\in [N]} \Delta_{h_1}\ldots \Delta_{h_k}f(x)
\]
\end{definition}

\begin{remark} One can define the Gowers norms on any abelian group; of special interest is the group $\F_2^n$ where the Gowers norms are intimately related to polynomial testing. 
\end{remark}

We make a few initial observations. For $1$-bounded functions $f$ ($\|f\|_{\infty} \le 1$)
\begin{itemize}
\item $\|f\|_{U^k[N]}=1$ if an only if   $f(x)=e^{2 \pi i q(x)}$ where $q$ is a polynomial of degree $<k$  . 
\item By repeated application of the Cauchy-Schwarz  inequality, if $f$ correlates with $e^{2 \pi i q(x)}$ where $q$ is a polynomial of degree $<k$ then $f$ has large Gowers norms; namely
\[
|\E_{x \in [N]} f(x)e^{-2 \pi i q(x)}|>\delta \implies
         \|f\|_{U^k[N]} \gg_{\delta} 1.
\]
\item If $f$ is a random function taking the values $\pm 1$ with probability $1/2$ for any $x\in [N]$, then by the law of large numbers,  $\|f\|_{U^k[N]}=o(1)$.
\end{itemize}
The Gowers uniformity norms play an important role in the study of arithmetic progressions. If $f$ and $g$ are close in the $U_k$ norm, i.e $\|f-g\|_{U^k[N]}$ is small, then they have approximately the same number of $k+1$ term progressions. Denote by $AP_{k}(f)$ the number of $(k+1)$-term progressions in $f$: denote
\[
AP_{k}(f)=\E_{x,d \in [N]} f(x)f(x+d)\ldots f(x+kd).
\]
Then 
 \begin{equation}\label{AP}
 |AP_{k}(f) -  AP_{k}(g)| \ll_k \|f-g\|_{U^k[N]}.
 \end{equation}
 \ \\
In fact a more general statement regarding linear forms is true:
\begin{proposition} Let $f_1, \ldots , f_k$ be $1$-bounded functions. Let $L_1(\vec x), \ldots, L_m(\vec x)$ be $k$  affine linear forms in $d$ variables with integer coefficients: $L_i(\vec x)=\sum_{j=1}^d l_{ij}x_j +b_i$, no two of which are affinely dependent. Then there exists $k>0$ such that 
\[
|\E_{\vec x \in [N]^d} f_1(L_1(\vec x)) \cdots f_k(L_m(\vec x)) | \ll \min_{j} \|f_j\|_{U^k[N]}.
\]
\end{proposition}
The proposition is proved via repeated applications of the Cauchy-Schwarz inequality, and this is  where the Gowers norms enter the picture in the proof of Theorem \ref{main}; it is the source of the condition that no two forms are affinely dependent. \\
 
The strategy of Gowers is similar in spirit to that of Roth. The idea is as follows. Let $E \subset [N]$ be with $|E|=\eta N$. Then 
\begin{itemize}
\item either the number of  {$(k+1)$}-term progressions is more than half of that expected in random set, namely $\ge \eta^{k+1}N^2/2$, or
\item $\|1_E -\eta\|_{U^{{k}}[N]} \gg_{\eta} 1$.
\end{itemize}

In order to proceed one needs to understand the condition $\|1_E -\eta \|_{U^{{k}}[N]} \gg_{\eta} 1$. For $k=2$ we observe that 
\[
\|f\|_{U^2[N]}^4 = \| \hat f \|_4^4 \le \|\hat f\|_2^2 \|\hat f\|_{\infty}^2.
\]
Thus if $\| f\|_2\le 1$ then we find that $\|f\|_{U^2[N]} \ge \eta$ implies $\|\hat f\|_{\infty} \ge \eta^2$. This implies that $f$ has a large Fourier coefficient, namely
\[
|\E_{x \in [N]} f(x)e(x\alpha)| \ge \eta^2.
\]
For larger $k$ the situation is much more complicated. Gowers proves the following local inverse theorem for higher Gowers norms.  

\begin{theorem}[Local inverse theorem for Gowers norms]\label{local} Let $f: \Z/N\Z \to \C$ be with $|f| \le 1$. Then
\[ \|f\|_{U^{k}[N]} \ge \delta \implies |\mathbb{E}_{x \in P} f(x)e^{2 \pi i q(x)}| \gg_{\delta} 1,\] where $P$ is a  progression of length at least $N^{t}$, $q(x)$ is a polynomial of degree ${k-1}$, and $t$ depends\footnote{In fact Gowers shows that one can find many such progressions: one can partition $\Z/N\Z$ into progressions $P_1, \ldots, P_M$ of average length greater than $N^{t}$, such that  $\sum_{i=1}^M|\sum_{x \in P_i} f(x)e^{2 \pi i q(x)}| \gg_{\delta} N$. } on $k, \delta$.
\end{theorem}

The word `local' in this context refers to the fact  that the  correlation in the above theorem is obtained not on the full interval $[N]$ but rather on a short progression of length at least $N^{t}$ with $t<1$ (for $k>2$).  This theorem provides sufficient structure to obtain increased density on a subprogression of length at least $N^{s}$: we apply 
Theorem \ref{local} to the function $1_E-\eta$, and use the equidistribution properties of the sequence $\{q(x)\}$  mod $1$ to find an arithmetic progression of length at least  $N^{s}$ ($s <t$) on which  $\{q(x)\}$mod $1$  is roughly constant.

\section{Classification of universal $k$-characteristic factors}
We return now to the question of classifying  $k$-characteristic factors. Recall that we are interested in  the averages 
\begin{equation}\label{multiple1}
\frac{1}{N}\sum_{n\le N} \int f(x)f(T^nx)f(T^{2n}x)\ldots  f(T^{kn}x) d\mu.
\end{equation}
The universal {$4$}-characteristic factors were classified by Host and Kra \cite{host-kra-personal}, and independently in the author's PhD thesis \cite{zieg-phd}, and were shown to be {$3$}-step pro-nilsystems. Both methods were extended to work for general $k$ - by Host and Kra in \cite{host-kra}, and by the author in \cite{zieg-jams}.
\begin{theorem}[Host-Kra (05), Z (07)]\label{nil-factor} Let ${\bf X}$ be an ergodic measure preserving system. 
The universal {$k$}-characteristic factor ${\bf Y_k}({\bf X})$  is a {$(k-1)$}-step pro-nilsystem.
\end{theorem}

We have the following diagram displaying the relation between the factors ${\bf Z_k}({\bf X})$ defined by Furstenberg in his proof of Szemer\'edi's theorem and the pro-nilfactors ${\bf Y_k}({\bf X})$ which are the universal characteristic factors:
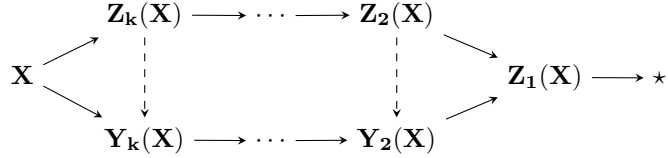
\begin{center}
\begin{tikzpicture}
  \matrix (m) [matrix of math nodes,row sep=0.7em,column sep=2em,minimum width=1em, ,ampersand replacement=\&]
  {
     {} \& {\bf Z_k}({\bf X}) \& {\cdots } \& {\bf Z_2}({\bf X}) \& {}  \& {}  \\
     {\bf X} \& {} \& {} \& {} \&  {\bf Z_1}({\bf X}) \& {\star}  \\
      {} \& {\bf Y_k}({\bf X}) \& { \cdots } \& {\bf Y_2}({\bf X}) \& {}  \& {}\\
     };
  \path[-stealth]
    (m-2-1) edge node [left] {} (m-1-2)
            edge  node [above] {} (m-3-2)
     (m-1-2) edge node [left] {} (m-1-3)
            edge  [dashed]   node [above] {} (m-3-2)      
      (m-1-3) edge node [left] {} (m-1-4)
      (m-3-3) edge node [left] {} (m-3-4)
      (m-3-4) edge node [left] {} (m-2-5)
      (m-3-2) edge node [left] {} (m-3-3)
            (m-2-5) edge node [left] {} (m-2-6)
    (m-1-4) edge node [right] {} (m-2-5)      
      edge  [dashed]   (m-3-4);
\end{tikzpicture}
\end{center}

As a corollary of this structure theorem one can calculate the asymptotic formula for the averages in \eqref{multiple1} via a limit formula for the corresponding averages on nilsystems \cite{zieg-limit}. 
\begin{theorem}[Z (05)]\label{limit}Let ${\bf X}$ be a $(k-1)$-step nilsystem. Then 
\[
\lim \frac{1}{N}\sum_{n\le N} \int f_0(x)f_1(T^nx) \ldots   f_{k}(T^{kn}x) d\mu = \int f_0(x_0)f_1(x_1) \ldots f_{k}(x_{k}) dm_H
\]
where $m_H$ is the Haar measure on the subnilmanifold $H\Gamma^{k+1}/\Gamma^{k+1} \subset X^{k+1}= (G/\Gamma)^{k+1}$, where $H$ is the subgroup 
\[
\{ (g_0, g_0g_1,g_0g_1^2g_2,g_0g_1^3g_2^3g_3, \ldots, g_0g_1^kg_2^{\binom{k}{2}} \ldots g_{k-2}^{\binom{k}{k-2}}) : g_i \in G_i \}
\]
where $\{1\}=G_{k-1} \subset G_{k-2} \subset \ldots \subset G_1 = G_0 = G$ is the derived series, i.e.  $G_{i+1}=[G_i,G]$. 
\end{theorem}

One can now prove Szemer\'edi's theorem by showing that the above limit is positive if $f_i=1_A$ for $i=0, \ldots, k$. This approach to proving Szemer\'edi's theorem (and various generalizations) was taken in \cite{bll}.\\

The proof in \cite{zieg-jams} generalizes the methods in \cite{cl1,cl2,cl3}. Inductively, one is led to the problem of solving a functional equation similar in nature to equation \eqref{cl}, only the extension cocycles are now defined on a (pro)-nilmanifolds (rather than a compact abelian group). Such cocycles are in general much more difficult to handle, but one can still use the fact that orbits on products of nilmanifolds are well understood and have a nice algebraic nature (as one can see in Theorem \ref{limit} above). \\

The proof in \cite{host-kra} introduces seminorms, which are similar, at least semantically, to the Gowers uniformity norms\footnote{Such averages as the one below were studied in the case $k=2$ already by Bergelson in \cite{Bergelson}.}
\begin{definition}[Host-Kra-Gowers semi-norms]
\[
\|f\|_{U^k({\bf X})}^{2^k}:= \lim_{N \to \infty } \mathbb{E}_{h_1, \ldots, h_k \in [N]} \int \Delta_{h_1} \ldots \Delta_{h_k} f(x) d\mu(x).
\]
\end{definition}
It is then proved that  characteristic factors for averages associated with the ergodic $U_k$ semi-norms defined above are also pro-nilsystems. Or, in a different formulation:
\[
\|f\|_{U^{k+1}({\bf X})} >0 \implies    \pi: {\bf X} \rightarrow \text{${k}$-step nilsystem},  \quad \pi_*f \ne 0.
\]
This suggests a far reaching generalization of the Gowers local inverse theorem, which we will discuss in section \ref{inverse} below. 

\section{The Green-Tao theorem from a characteristic factor point of view}
In their famous paper, Green and Tao prove a Szemer\'edi's type theorem in the prime numbers \cite{gt-primes}: 
\begin{theorem}[Green-Tao (05)]\label{gt-primes}
Let $E \subset \mathbb P$ of positive relative density, then $E$ contains long arithmetic progressions.
\end{theorem}

We present the idea of the proof from a characteristic factor point of view. Our starting point will be the following version of Szemer\'edi's theorem: Let $f:[N] \to [0,1]$ be a function with $|\E_{n \in [N]} f(x)| >\delta$. Then 
for any integer $k>0$
\begin{equation}\label{comb-average}
AP_{k}(f)=\E_{x,d \in [N]} f(x)f(x+d)\ldots f(x+kd) \ge c(\delta)+o(1)
\end{equation}
where $c(\delta)>0$, and is independent of $N$. 

If we naively try to apply this theorem for a subset $E$ of the prime numbers of relative density $\delta$,  we run into an obvious problem that $\E_{x \in [N]} 1_E(x)=o(1)$. We can try to fix this problem by putting a weight on each prime - we consider the von-Mangoldt function $\Lambda(x)$ which takes the value $\log p$ if $x$ is a positive power of $p$ and $0$ otherwise. In this case we will have 
\[
\E_{x \in [N]} \Lambda(x) 1_E(x)=\delta +o(1),
\]
for some constant $\delta>0$ (independent of $N$).
But now we face the problem that the function $\tilde 1_E(x)= \Lambda(x) 1_E(x)$ does not take values in $[0,1]$; in fact the function $\tilde 1_E(x)$ is unbounded. Green and Tao show that for a certain class of unbounded functions (functions bounded by a $k$-psuedorandom function) one can find a "$k$-characteristic factor"  for the average \eqref{comb-average} generated by {\em bounded} functions !   We can summarize the procedure as follows:
\begin{itemize}
\item  Introduce combinatorial notions of (approximate) {\em factor} and {\em  projection} onto a factor. 
\item  Find a convenient combinatorial {"$k$-characteristic factor"} for averages associated with the $U_k$ norms, in this case a factor of functions bounded by a constant $C(k)$ (depending only on $k$).
\item    Let $\pi_* (\tilde1_E )$ be  the (approximate) projection on the factor.  Then $0 \le \pi_* (\tilde1_E ) \le C(k)$, the average of the function $\pi_* \tilde1_E $ is approximately the same as that of $\tilde1_E$, namely approximately $\delta$, and  $\| \tilde1_E-\pi_* (\tilde1_E )\|_{U^k[N]}$\footnote{We defined the Gowers norms for functions $f$ on the group $\Z/N\Z$. We can define Gowers norms for functions $f : [N] \rightarrow \C$, setting $G := \Z/\tilde N\Z$ for some integer $\tilde N \geq 2^d N$, and defining  a function $\tilde f : G \rightarrow \C$ by $\tilde f(x) = f(x)$ for $x = 1,\dots,N$ and $\tilde f(x) = 0$ otherwise. We then set 
\[ \Vert f \Vert_{U^d[N]} := \Vert \tilde f \Vert_{U^d(G)} / \Vert 1_{[N]} \Vert_{U^d(G)},\] where $1_{[N]}$ is the indicator function of $[N]$.  It is easy to see that this definition is independent of the choice of $\tilde N$.} is small. A version of the Gowers-Cauchy-Schwarz  inequality (for functions bounded by $k$-pseudoradnom functions) gives then, as in \eqref{AP}, that
 \begin{equation}\label{AP-prime}
 |AP_{k}(\tilde1_E) -  AP_{k}(\pi_*\tilde1_E)| \ll \|\tilde1_E- \pi_* \tilde1_E\|_{U^k[N]}
 \end{equation}
\item Apply Szemer\'edi's Theorem to the $C(k)$-bounded function $\pi_* \tilde 1_E $, to obtain $AP_{k}(\pi_*\tilde1_E) \gg_{\delta} 1$, and thus 
$AP_{k}(\tilde1_E) \gg_{\delta} 1$
\end{itemize}
\ \\
A different way to say this is that given $\epsilon>0$ we can decompose 
\begin{equation}\label{decomposition}
\tilde1_E = g+h
\end{equation}
where $g$ is a $C(k)$-bounded function, and $h$ is a function with $\|h\|_{U^k[N]}< \epsilon$. This type of theorem is now referred to as a {\em decomposition} theorem. There is a very nice modern and more abstract treatment of general decomposition theorems in \cite{gowers-reg}, and \cite{rtts} using the Hahn-Banach theorem.  We remark that Theorem \ref{gt-primes} has since been extended to include polynomial configurations \cite{tao-ziegler-polyprimes}, and multidimensional configurations \cite{tao-ziegler-multi, fz,cmt}.

\section{The Inverse Theorem for the Gowers Norms}\label{inverse}

The argument in the Green-Tao theorem is based on Szemer\'edi's theorem which is valid for {\em any} subset of positive density in the integers. This has two major caveats. The first is that it can not lead to an asymptotic formula for the number of arithmetic progressions, only a lower bound. The second is that it can not be used to study non homogeneous linear configurations, since there are counter examples within periodic sets of positive density. How then can we hope to get an asymptotic formula as in Theorem \ref{main} ?\\

We recall now that - in the ergodic theoretic context - to get a limit formula we needed to identify the universal characteristic factors.    
Motivated by theorem \ref{nil-factor},  Green and Tao conjectured in 2006 that the combinatorial ``universal characteristic factors" for the $U_k$ norm  
come from {\em nilsequences} - sequences arising in a natural way from nilsystems. 

\begin{conjecture}[Inverse conjecture for the Gowers norms ($\GI(s)$)]\label{conj-main}  Let $s \geq 0$ be an integer, and let $0 < \delta \leq 1$.   Then there exists a finite collection ${\mathcal M}_{s,\delta}$ of $s$-step nilmanifolds $G/\Gamma$, each equipped with some smooth Riemannian metric $d_{G/\Gamma}$ as well as constants $C(s,\delta), c(s,\delta) > 0$ with the following property. Whenever $N \geq 1$ and $f : [N] \rightarrow \C$ is a $1$-bounded function  such that $\Vert f \Vert_{U^{s+1}[N]} \geq \delta$, there exists a nilmanifold $G/\Gamma \in {\mathcal M}_{s,\delta}$, some $g \in G$ and a function $F: G/\Gamma \to \C$ bounded in magnitude by $1$ and with Lipschitz constant at most $C(s,\delta)$ with respect to the metric $d_{G/\Gamma}$, such that
\begin{equation}\label{correlation-nil} |\E_{n \in [N]} f(n) \overline{F(g^n x)}| \geq c(s,\delta).
\end{equation}
\end{conjecture}

That is, the global obstruction (scale $N$) to Gowers uniformity come from sequences arising from nilsystems. 
Recall that the local theorem for the Gowers norms shows that
local obstructions (at scale $N^{t}$) to Gowers $U^{s+1}$ uniformity norms come from phase polynomials of degree ${s}$. We remark that the converse to Conjecture \ref{conj-main} is true and relatively easy to prove via repeated applications of the Cauchy-Schwarz inequality. Namely, if \eqref{correlation-nil} holds then $\|f\|_{U^{s+1}[N]} \gg_{\delta} 1$.  We also mention that if $\delta$ is sufficiently close to $1$ then the conjecture is true; moreover,  $f$ is close (in $L^1$) to a genuine (unique) phase polynomial \cite{akklr}, and thus correlates with a (unique) phase polynomial\footnote{One can exhibit a polynomial phase function $e^{p(x)}$ as a nilsequence see e.g. \cite{gtz}}; uniqueness allows one to try an intelligent guess. In the realm when $\delta>0$, we cannot expect uniqueness, and as it turns out, we also can't expect correlation with a genuine phase polynomial.

One can ask a similar question in the context of finite field geometry. Given a function $f:\F_p^n \to \D$ with large Gowers norm (fixing $p$ and letting $n$ approach $\infty$), what can be said about $f$? It was conjectured  that such functions would correlate with polynomial phase functions. More precisely:
\begin{conjecture}[Inverse conjecture for the Gowers norms in finite fields]\label{conj-ff}
Let $p$ be a prime and let $f:\F_p^n \to \C$ be $1$-bounded, with $\|f\|_{U^{s+1}[\F_p^n]}\ge  \delta$. Then there exists a polynomial $P:\F_p^n \to \F_p$ of degree $\le k$ such that  
\[
 |\mathbb{E}_{x\in \F_p^n}f(x)e^{2\pi i P(x)/p}| \ge c(s,\delta).
\]\end{conjecture}

The case $s=1$ of both conjectures follows from a short Fourier-analytic computation.
The case $s=2$  of Conjecture \ref{conj-main} was proved in \cite{gt:inverse-u3}.  The case $s=2$  of Conjecture \ref{conj-ff} was proved in \cite{gt:inverse-u3} for odd $p$ and for $p=2$  in \cite{sam}. Surprisingly, Conjecture \ref{conj-ff} turned out to be {\em false}; a counter example for the $U^4[\F_2^n]$ was constructed independently in  \cite{gt-false,lms}. However, it turned out that with a small modification ofConjecture \ref{conj-ff} is actually true \cite{fp-inverse,tao-ziegler-correspondence,tao-ziegler-lowchar}.
Call $P:\F_p^n \to \C$ a non-standard polynomial of degree $<k$ if for all $h_1, \ldots, h_s \in \F_p^n$ we have 
\[
\Delta_{h_1}\ldots \Delta_{h_s} P \equiv 1
\]
If char $\F\ge s$, then a non-standard polynomial is a standard phase polynomial, i.e $e^{2\pi i P(x)/p}$ where $P:\F_p^n \to \F_p$ a polynomial of degree $<s$, but otherwise the class of non-standard polynomials is larger.

\begin{theorem}[Bergelson-Tao-Z (10), Tao-Z (10,12)]\label{btz}
Let $p$ be a prime and let $f:\F_p^n \to \mathbb C$ be $1$-bounded, with $\|f\|_{U^{s+1}[\F_p^n]}\ge \delta$. Then there exists a non-standard polynomial $P$ of degree $\le s$, and a constant  $c(s, \delta)>0$ such that 
\[
|\mathbb{E}_{x\le \F_p^n}f(x)e^{P(x)}| \gg c(s, \delta).
\]
\end{theorem}
The proof of theorem \ref{btz} is via an ergodic theoretic structure theorem, similar in nature to Theorem  \ref{nil-factor}, and a correspondence theorem - translating the finitary question to a question about limiting behavior of multiple averages for an $\oplus \F_p$ ergodic action. 

Finally Conjecture \ref{conj-main} was proved \cite{gtz}:
\begin{theorem}[Green-Tao-Z (12)]\label{gtz}
The inverse conjecture for the Gowers norms $GI(s)$ norms is true.
\end{theorem}

The proof of Theorem  \ref{gtz} is long an complicated and is carried out in \cite{gtz}. For a more gentle introduction to the proof we refer the reader to either \cite{gtz2}, where the case $k=3$ (the $U^4$ norm) is handled,  or to the announcement in \cite{gtz1}. We now try to give the flavor of the proof.  Suppose $\Vert f \Vert_{U^{s+1}[N]} \geq \delta$, then by definition
\[
\E_{h \in N} \|\Delta_h f(n)\|^{2^s}_{U^{s}[N]}  \gg_{\delta} 1.
\]
It follows that for all $h$ in a set $H$ of size  $\gg_{\delta} N$ we have  $\|\Delta_h f(n)\|_{U^{s}[N]}  \gg_{\delta} 1$.  Now, inductively we know that,
for $h \in H$, $\Delta_h f(n)$ correlates with an $s-1$-step nilsequence $F_h(g_h^nx_h\Gamma)$ (of complexity $\ll_{\delta } 1$), namely
\[
|\E_{h \in N} \Delta_h f(n)F_h(g_h^nx_h\Gamma)| \gg_{\delta} 1
\] 
In the case $GI(2)$, this $1$-step nilsequence can be taken to be $e^{2\pi i \lambda_h n}$, but in general we can't hope for anything as simple. The key difficulty now is to try to find some extra structure relating the nilsequences $F_h(g_h^nx_h\Gamma)$ for different values of $h$. This is already quite difficult in the $GI(2)$ case. In this case, an ingenious argument of Gowers involving tools from additive combinatorics, coupled with some geometry of numbers allows one to linearize $\lambda_h$ on a nice set - a generalized arithmetic progression (GAP). This argument is then combined with a symmetry argument to construct a $2$-step nilsequence $g(h)$ with  $\Delta_hg(n)=e^{2\pi i \lambda_h n}$ for many values of $h$ \cite{gt:inverse-u3}). For general $s$, we follow the same strategy, however it turns out to be rather difficult to extract some algebraic structure relating the various nilsequences $F_h(g_h^nx_h\Gamma)$.  An alternate approach to the inverse theorem was subsequently developed by Szegedy  \cite{ camarena-szegedy, szegedy}. We remark that both Theorems \ref{btz}, \ref{gtz} are qualitative; it is a major open question to find quantitative proofs for them.\\

How can one apply Theorem \ref{gtz} to obtain Theorem \ref{main}? We give a very rough sketch. One needs to calculate the projection of the function $\tilde 1_{\mathbb{P}}(n) =  (\log n) 1_{\mathbb{P}}(n)$  onto the combinatorial nil-factor.  It turns out that the projection essentially lies in the much smaller factor of periodic functions (with bounded period). More precisely, one first performs pre-sieving to eliminate the periodic contributions. Let $W=\prod_{p<w} p$ for $w$ a slowly increasing function of $N$.  For $(b,W)=1$ consider $\tilde 1_{W,b,\mathbb{P}}(n) = \tilde 1_{\mathbb{P}}(Wn+b) $.  The projection of this function on the combinatorial nil-factor should be constant, and since its average is $1$ - this constant should be $1$; namely one must show that
\[
\|\tilde 1_{W,b,\mathbb{P}}(n)  - 1\|_{U^k[N]} =o(1).
\]
Suppose $\|\tilde 1_{W,b,\mathbb{P}}(n)  - 1\|_{U^k[N]}>\delta$. Fix $\eps>0$, and decompose as in \eqref{decomposition}
\[
 \tilde 1_{W,b,\mathbb{P}}(x)  - 1= f+g
\]
where $f$ is bounded $f \ll_k 1$ , and $\|g\|_{U^k[N]}< \eps$. Then since $U_k$ is a norm we get that
$\|f\|_{U^k[N]}>\delta/2$. By Theorem \ref{gtz} there is a nilsequence $F(g^nx\Gamma)$ of complexity $\ll_{\delta} 1$ (i.e. all parameters associated  with the nilsequence such as the dimension of the nilmanifold are bounded in terms of $\delta$),  such that  $|\E f(x)F(g^nx\Gamma)| \gg_{\delta} 1$.
From the easy direction of Theorem \ref{gtz} (which is valid for non bounded functions as well, via repeated applications of the Cauchy-Schwarz inequality) we have $|\E g(x)F(g^nx\Gamma)| < c(\eps)$ (with $c$ a decreasing function). In \cite{gt-mobius} it is shown 
that for any bounded complexity nilsequence we have  $\E (\tilde 1_{W,b,\mathbb{P}}(n)  - 1 )F(g^nx\Gamma)=o(1)$. Choosing $\eps$ sufficiently small in the decomposition \eqref{decomposition}, we get a contradiction. 

\section{Acknowledgement}

I thank H. Furstenberg for introducing me to ergodic theory and to the rich subject of multiple recurrence. 
I  thank V. Bergelson and T. Tao for their valuable comments on an earlier version of this paper.

\end{document}